\documentclass[a4paper,10pt,conference]{ieeeconf}
\IEEEoverridecommandlockouts
\usepackage[utf8x]{inputenc}
\DeclareUnicodeCharacter{"FF0E}{}
\usepackage[dvips]{graphicx}
\usepackage{mac_ECC2021}
\usepackage{amsmath}
\usepackage{amssymb}
\usepackage{fleqn}
\usepackage{yhmath}
\usepackage{noitemsep}
\usepackage{here}
\usepackage{pstricks}
\usepackage{footnote}
\usepackage{mathfonts}
\usepackage{psfrag}
\usepackage{arydshln}
\usepackage{macros_e}

\everymath={\displaystyle}
\title{
\LARGE \bf 
$l_2$ Induced Norm Analysis
of Discrete-Time LTI Systems for
Nonnegative Input Signals  
and Its Application to \\
Stability Analysis of Recurrent Neural Networks
}
\author{Yoshio Ebihara, Hayato Waki, Victor Magron, \\
Ngoc Hoang Anh Mai, Dimitri Peaucelle, Sophie Tarbouriech
\thanks{
Y. Ebihara is with the 
Graduate School of Information Science and 
Electrical Engineering, Kyushu University, 
744 Motooka, Nishi-ku, Fukuoka 819-0395, Japan, 
he was also with 
LAAS-CNRS, Universit\'{e} de Toulouse, CNRS, Toulouse, France,   
in 2011.  
H. Waki is with the Institute of Mathematics for Industry,
Kyushu University, 744 Motooka, Nishi-ku, Fukuoka 819-0395, Japan.  
V. Magron, N. H. Mai, D. Peaucelle, and S. Tarbouriech are
LAAS-CNRS, Universit\'{e} de Toulouse, CNRS, Toulouse, France.  
}%
}

\begin{document}
\maketitle
\thispagestyle{empty}
\pagestyle{empty}

\begin{abstract}
In this paper, we focus on the ``positive'' $l_2$ induced norm of 
discrete-time linear time-invariant systems where the input signals are
restricted to be nonnegative.  
To cope with the nonnegativity of the input signals, 
we employ copositive programming as the mathematical tool for the
analysis.  
Then, by applying an inner approximation to the copositive cone, 
we derive numerically tractable 
semidefinite programming problems for the upper and lower bound
computation of the ``positive'' $l_2$ induced norm.  
This norm is typically useful for the stability analysis of
feedback systems constructed from an LTI system and 
nonlinearities where the nonlinear elements
provide only nonnegative signals.  
As a concrete example, 
we illustrate the usefulness of the ``positive'' $l_2$ induced norm
for the stability analysis of recurrent neural networks
with activation functions being rectified linear units.  

\end{abstract}


\section{Introduction}

The $l_2$ ($L_2$) induced norm plays a central role
in stability and performance analysis of discrete-time (continuous-time)
feedback systems \cite{Khalil_2002}.  
As is well known, 
small-gain stability criterion allows us to 
assess the stability of feedback systems
constructed from two subsystems by evaluating
their $l_2$ induced norms.  
One of the key discoveries on the $l_2$ induced norm 
of (finite-dimensional) linear time-invariant (LTI) systems
would be KYP lemma \cite{Rantzer_SCL1996}, which
characterizes the induced norm 
by a semidefinite programming problem (SDP).   
It should be noted that, 
even though the $l_2$ induce norm is defined in time-domain, 
the core in deriving KYP lemma is the
treatments of LTI systems in frequency-domain.  

In the standard $l_2$ induced norm analysis of 
LTI systems, we of course presume that 
the input signals are sign indefinite in time-domain.  
On the contrary, in this paper, we focus 
on the ``positive'' $l_2$ induced norm of 
discrete-time LTI systems where the input signals are
restricted to be nonnegative.  
As clarified in this paper, 
this norm is typically useful for the stability analysis of
feedback systems constructed from an LTI system and 
nonlinearities where nonlinear elements
provide only nonnegative signals.  
It is also true that the norm analysis is 
partly motivated from our preceding studies 
on positive systems \cite{Ebihara_IEEE2017, Kato_LCSS2020}, 
where the treatments of nonnegative signals are essentially important.  

The analysis of the ``positive'' $l_2$ induced norm 
for discrete-time LTI systems is mathematically challenging 
due to the following reasons:
\begin{itemize}
 \item[(i)] The nonnegativity constraint on the input signals 
	     is a genuine time-domain constraint and hence it does not
	     allow us to carry out the analysis in frequency-domain.  
 \item[(ii)] Even though SDP is commonly used 
	    for LTI system analysis, 
	    the positive semidefinite cone employed in SDP
	    has no functionality to distinguish nonnegative vectors or signals.  
\end{itemize}
To get around these difficulties and 
cope with the nonnegativity of the input signals, 
we loosen the positive semidefinite cone 
to the copositive cone 
and employ copositive programming (COP) \cite{Duer_2010}  
as a mathematical tool for the analysis.  
COP is a convex optimization problem 
on the positive semidefinite cone, 
but unfortunately known to be numerically intractable \cite{Duer_2010}.   
Therefore, by further applying an inner approximation to the copositive cone, 
we derive numerically tractable 
SDPs for upper and lower bound
computation of the ``positive'' $l_2$ induced norm.  
We illustrate the usefulness of the ``positive'' $l_2$ induced norm
for the stability analysis of recurrent neural networks (RNNs)
with activation functions being rectified linear units (ReLUs).  
Recently, the usefulness of 
RNNs is widely recognized for analysis and estimation
of time series generated by (hidden) dynamical systems.   
This is achieved by incorporating feedback loops in the networks. 
However, the existence of feedback loops could be a source of
instability, and stability analysis of RNNs remains to be an 
outstanding issue in the fields of neural network and machine learning
\cite{Barabanov_IEEENN2002,Zhang_IEEENN2014,Salehinejad_2018}.  
By making use of the fact that ReLUs only provide nonnegative signals, 
we derive novel small-gain type stability conditions for RNNs
on the basis of the ``positive'' $l_2$ induced norm.  
We illustrate the effectiveness of the new stability tests
by numerical examples.

We use the following notation in this paper.
The set of natural numbers is denoted by $\bbN$.  
The set of $n\times m$ real matrices is denoted by $\bbR^{n\times m}$, and
the set of $n\times m$ entrywise nonnegative (strictly positive) matrices is denoted
by $\bbR_+^{n\times m}\ (\bbR_{++}^{n\times m})$.
For a matrix $A$, we also write $A\geq 0\ (A>0)$ to denote that  
$A$ is entrywise nonnegative (strictly positive).   
For $A\in\bbR^{n\times n}$, we define $\He\{A\}:=A+A^T$.

We denote the set of $n\times n$ real symmetric matrices by $\bbS_n$, and
define the inner product of $A, B\in\bbS_n$ by 
$\langle A,B\rangle :=\mathrm{trace}(AB)$.
For a set $\clK_n\subset\bbS_n$, the interior of $\clK_n$ is denoted by 
$\clK_n^\circ$.
The dual cone of $\clK_n\subset\bbS_n$ is defined by
$\mathcal{K}_n^\ast:=\{A\in\bbS_n|\ \forall B\in\mathcal{K}_n,\ 
	\langle A,B\rangle\geq 0\}$.  
If $\clK_n\subset\bbS_n$ is a closed convex cone 
which has nonempty interior and pointed 
(i.e., $x,-x\in\clK_n\ \Ra\ x=0$), 
then $\clK_n$ is called proper \cite{Boyd_2004}.
The interior of a proper cone $\clK_n$ is defined as follows
\cite{Berman_1979}:  
\begin{align}\label{coneint}
  \mathcal{K}_n^\circ=\{A\in\bbS_n|\ \forall
	B\in\mathcal{K}_n^\ast\backslash\{0\},\ \langle A,B\rangle >0\}.   
\end{align}
For a proper cone $\clK_n\subset\bbS_n$ and $P\in\bbS_n$, 
by following \cite{Boyd_2004}, 
we define a partial ordering $P\succeq_\clK 0\ (P\preceq_\clK 0)$ which means
$P\in\clK_n\ (-P\in\clK_n)$.
We also write $P\succ_\clK 0\ (P\prec_\clK 0)$ if 
$P\in\clK_n^\circ\ (-P\in\clK_n^\circ)$.
The subscript $\clK$ is often omitted if $\clK$ is 
the positive semidefinite cone.  
For an $\bbS_n$-valued affine function $F(x)$, we call the inequalities
of the form 
$F(x)\succeq_\clK 0\ (F(x)\preceq_\clK 0)$ and $F(x)\succ_\clK 0\
(F(x)\prec_\clK 0)$ the 
linear matrix inequalities (LMIs) on the cone $\clK$.

For a discrete-time signal $w$ defined 
over the time interval $[0,\infty)$, we define
\[
\begin{array}{@{}l}
 \|w\|_{2}:=\sqrt{\sum_{k=0}^{\infty}|w(k)|_2^2}
\end{array}
\]
where for $v\in\bbR^{n_v}$ we define
$|v|_2:=\sqrt{\sum_{j=1}^{n_v} v_j^2}$.  
We also define
\[
\begin{array}{@{}l}
 l_{2}  :=\left\{w:\ \|w\|_{2}<\infty \right\},\\ 
 l_{2+} :=\left\{w:\ w\in l_{2},\ w(k)\ge 0\ (\forall k\ge 0)\right\}  
\end{array}
\]
and
\[
\begin{array}{@{}l}
 l_{2e}  :=\left\{w:\ w_\tau\in l_2,\ \forall \tau\in[0,\infty) \right\}
\end{array}
\]
where $w_\tau$ is the truncation of the signal $w$ up to the time instant
$\tau$
and defined by
\[
 w_\tau(k)=
 \left\{
 \begin{array}{cc}
  w(k)& (k\le \tau), \\
  0& (k> \tau).  \\
 \end{array}
 \right.  
\]
For an operator $H:\ l_{2e}\ni w \to z \in l_{2e}$, 
we define its $l_2$ induced norm by
\begin{equation}
\|H\|_{2}:=\sup_{w\in l_2,\ \|w\|_2=1} \ \|z\|_2.  
\label{eq:l2norm}
\end{equation}
The main interest of the present paper concerns
a variant of the $l_2$ induced norm 
with input being restricted to be nonnegative which is defined by
\begin{equation}
\|H\|_{2+}:=\sup_{w\in l_{2+},\ \|w\|_2=1} \ \|z\|_2.  
\label{eq:l2+norm}
\end{equation}
With a little abuse of notation, for a matrix $M\in\bbR^{n\times m}$, 
we define
\begin{equation}
\begin{array}{@{}l}
\|M\|_{2}:=\sup_{v\in \bbR^m,\ |v|_2=1} \ |Mv|_2\quad (= \smax(M)),\\
\|M\|_{2+}:=\sup_{v\in \bbR_+^m,\ |v|_2=1} \ |Mv|_2  
\end{array}
\label{eq:possing}
\end{equation}
where $\smax(M)$ is the maximal singular value of $M$.

\section{Copositive Programming Problem (COP)}\label{sec:cop}

Copositive Programming Problem (COP) is a convex optimization problem in which 
we minimize a linear
objective function over the LMI constraints on the copositive cone   
\cite{Duer_2010}.  
In this section, we summarize its basics.  

\subsection{Convex Cones Related to COP}

Let us review the definition and the property of convex
cones related to COP.
\begin{definition}\cite{Berman_2003}\label{conedef} 
The definition of proper cones 
$\PSD_n$, $\COP_n$, $\CP_n$, $\NN_n$, and $\DNN_n$ in 
$\bbS_n$ are as follows.
\begin{enumerate}
  \item 
  $\PSD_n:=\{P\in\bbS_n:\ \forall x\in\bbR^n,\
	x^TPx\geq 0\}=\{P\in\bbS_n:\ \exists B\
	\mbox{s.t.}\ P=BB^T\}$ is called {\it the positive semidefinite cone}.
  \item
  $\COP_n:=\{P\in\bbS_n:\ \forall x\in\bbR_{+}^n,\
	x^TPx\geq 0\}$ is called {\it the copositive cone}.
  \item
  $\CP_n:=\{P\in\bbS_n:\ \exists B\ge 0\
	\mbox{s.t.}\ P=BB^T\}$ is called {\it the completely positive cone}.
  \item
  $\NN_n:=
  \{P\in\bbS_n:\ P\geq 0\}
  $
  is called {\it the nonnegative cone}.
  \item
  $\PSD_n+\NN_n:=\{P+Q:\ P\in\PSD_n,\
	Q\in\NN_n\}$．
       This is the Minkowski sum of the positive semidefinite cone and 
       the nonnegative cone.
  \item
  $\DNN_n:=\PSD_n\cap\NN_n$ is called 
  {\it the doubly nonnegative cone}.
\end{enumerate}
\end{definition}

From Definition \ref{conedef}, 
we clearly see that the following inclusion relationships hold: 
\[
 \CP_n\subset\DNN_n\subset\PSD_n\subset\PSD_n+\NN_n\subset\COP_n,
\]
\[
 \CP_n\subset\DNN_n\subset\NN_n\subset\PSD_n+\NN_n\subset\COP_n.
\]
In particular, when $n\leq 4$, it is known that 
$\COP_n=\PSD_n
+\NN_n$ and $\CP_n =\DNN_n$ hold \cite{Berman_2003}．
On the other hand, on the duality of these cones, 
$\COP_n$ and $\CP_n$ are dual to each other, 
$\PSD_n+\NN_n$ and $\DNN_n$ are dual to each other, 
and $\PSD_n$ and $\NN_n$ are
self-dual.

When dealing with analysis problems of LTI systems, 
we often need to consider strict LMI conditions.  
With this fact in mind, let us review 
the characterization of 
the interiors of the convex cones in Definition~\ref{conedef}.  
\begin{proposition}\cite{Duer_EJLA2008,Dickinson_EJLA2010}\label{copint}
The interiors of the proper cones given in Definition \ref{conedef} are
characterized as follows.
\begin{enumerate}
  \item
  $\PSD_n^\circ=\{P\in\bbS_n:\  \forall
	x\in\bbR^n\backslash\{0\},\ x^TPx>0\}=
  \{P\in\bbS_n:\  \exists B\ \mbox{s.t.}\ P=BB^T,\
       \mathrm{rank}(B)=n\}$.  
  \item
  $\COP_n^\circ=\{P\in\bbS_n:\  \forall
	x\in\bbR_{+}^n\backslash\{0\},\ x^TPx>0\}$.  
  \item
  $\CP_n^\circ=\{P\in\bbS_n:\  
	\exists B\ \mbox{s.t.}\ P=BB^T,\ B>0,\ \mathrm{rank}(B)=n\}$.  
  \item
  $\NN_n^\circ=  \{P\in\bbS_n:\  P>0\}.  
  $
  \item
  $(\PSD_n+\NN_n)^\circ
  =\PSD_n^\circ+\NN_n^\circ  
  $.  
  \item
  $\DNN_n^\circ=\PSD_n^\circ\cap\NN_n^\circ$.  
\end{enumerate}
\end{proposition}

\subsection{Basic Properties of COP}

COP is a convex optimization problem on the copositive cone and
its dual is a convex optimization problem on the completely positive cone.  
As mentioned in \cite{Duer_2010}, 
the problem to determine whether a given symmetric matrix is
copositive or not is a co-NP complete problem, and
the problem to determine whether a given symmetric matrix is
completely positive or not is an NP-hard problem.  
Therefore, it is hard to solve COP 
numerically in general.  
However, since the problem to determine whether a given matrix is in $\PSD+\NN$
or in $\DNN$ can readily be reduced to SDPs, 
we can numerically solve the convex optimization problems on the cones
$\PSD+\NN$ and $\DNN$ easily.  
Moreover, when $n\leq 4$, it is known that $\COP_n=\PSD_n+\NN_n$ and 
$\CP_n=\DNN_n$ as stated above, 
and hence those COPs with $n\leq 4$ can be reduced to SDPs.  

\section{$l_2$ Induced Norm Analysis of LTI Systems for Nonnegative Inputs}
\label{sec:l2}

\subsection{Problem Description}
\label{sub:Problem}
Let us consider the discrete-time LTI system $G$ given by
\begin{equation}
G:\ 
\left\{
\arraycolsep=0.5mm
\begin{array}{cccccccc}
 x(k+1)&=& A x(k)& + &B w(k), &\ x(0)=0,\\
   z(k)&=& C x(k)& + &D w(k) \\
\end{array}
\right.  
\label{eq:G}
\end{equation}
where
$A\in\bbR^{n\times n}$, 
$B\in\bbR^{n\times \nw}$,
$C\in\bbR^{\nz\times n}$, and
$D\in\bbR^{\nz\times \nw}$.   
We assume that the system $G$ is stable, i.e., 
the matrix $A$ is Schur-Cohn stable.  
It is well known that the $l_2$ induced norm 
$\|G\|_2$ defined by \rec{eq:l2norm} coincides with the
$H_\infty$ norm for stable LTI systems
and plays an essential role in stability analysis of feedback systems.  
In this paper, we are interested in computing 
the $l_2$ induced norm where the input signal $w$ is 
constrained to be nonnegative.  
Namely, we focus on the computation of 
the ``positive'' $l_2$ induced norm 
$\|G\|_{2+}$ defined by
\rec{eq:l2+norm}.  
From the definition \rec{eq:l2+norm}, it is very clear that
$\|G\|_{2+}\le \|G\|_{2}$.  
Here, note that 
a discrete-time LTI system of the form \rec{eq:G}
is said to be externally positive if its output is 
nonnegative for any nonnegative input under zero initial state 
\cite{Farina_2000}.  
Then, in this case, it is well known that 
$\|G\|_{2+}= \|G\|_{2}$, see, e.g.,  \cite{Tanaka_IEEE2011,Rantzer_IEEE2016}.  

\subsection{Basic Results}
\label{sub:basic}

The next result forms an important basis of this study.  
\begin{theorem}
For the stable LTI system $G$ described by \rec{eq:G}
and given $\gamma>0$, suppose 
there exist $P\in\PSD_n$ and $Q\in\COP_{\nw}$
such that
\begin{equation}
\begin{array}{@{}l}
L(A,B,C,D,P,Q,\gamma)\prec 0  
\end{array}
\label{eq:LMI}
\end{equation}
where
\[
\begin{array}{@{}l}
L(A,B,C,D,P,Q,\gamma)\\:= \left[
\arraycolsep=0.5mm
\begin{array}{cc}
 -P & 0 \\
 0 & -\gamma^2 I_{\nw}+Q
\end{array}
\right]+
 \left[
 \begin{array}{cc}
  A & B \\
  C & D\\
 \end{array}
 \right]^T
 \left[
 \begin{array}{cc}
  P & 0\\
  0 & I_{\nz}\\
 \end{array}
 \right]
 \left[
 \begin{array}{cc}
  A & B \\
  C & D\\
 \end{array}
 \right].  
\end{array}
\]
Then we have $\|G\|_{2+}<\gamma$.  
\label{th:main}
\end{theorem}
\begin{proofof}{\rth{th:main}}
We first note that if \rec{eq:LMI} holds 
with $P$ and $Q$ then
there exists a sufficiently small $\varepsilon>0$ such that
the next condition holds with exactly the same $P$, $Q$
and $\tilgam:=\gamma-\varepsilon$:
\begin{equation}
\scalebox{1.0}{$
\begin{array}{@{}l}
\arraycolsep=0.5mm
\left[
\begin{array}{cc}
 -P & 0 \\
 0 & -\tilgam^2 I_{\nw}+Q
\end{array}
\right]+
 \left[
 \begin{array}{cc}
  A & B \\
  C & D\\
 \end{array}
 \right]^T
 \left[
 \begin{array}{cc}
  P & 0\\
  0 & I_{\nz}\\
 \end{array}
 \right]
 \left[
 \begin{array}{cc}
  A & B \\
  C & D\\
 \end{array}
 \right]\prec 0.  
\end{array}$}
\label{eq:LMI2}
\end{equation}
With this fact in mind, 
let us consider the trajectory of the state $x$
corresponding to the input signal $w\in l_{2+}$
with $\|w\|_2=1$ for the system $G$.  
Then, from \rec{eq:LMI2}, we have
\[
\begin{array}{@{}l}
\scalebox{0.76}{$
\begin{array}{@{}l}
\arraycolsep=0.3mm
\left[
\begin{array}{c}
x(k)\\
w(k)\\
\end{array}
\right]^T\left\{
\left[
\begin{array}{cc}
 -P & 0 \\
 0 & -\tilgam^2 I_{\nw}+Q
\end{array}
\right]+
 \left[
 \begin{array}{cc}
  A & B \\
  C & D\\
 \end{array}
 \right]^T
 \left[
 \begin{array}{cc}
  P & 0\\
  0 & I_{\nz}\\
 \end{array}
 \right]
 \left[
 \begin{array}{cc}
  A & B \\
  C & D\\
 \end{array}
 \right]\right\}
\left[
\begin{array}{c}
x(k)\\
w(k)\\
\end{array}
\right]\le 0
\end{array}$},\\
(k=0,1,\cdots)
\end{array}
\]
or equivalently, 
\[
\begin{array}{@{}l}
\scalebox{0.72}{$
\begin{array}{@{}l}
\arraycolsep=0.3mm
\left[
\begin{array}{c}
x(k)\\
w(k)\\
\end{array}
\right]^T
\left[
\begin{array}{cc}
 -P & 0 \\
 0 & -\tilgam^2 I_{\nw}+Q
\end{array}
\right]
\left[
\begin{array}{c}
x(k)\\
w(k)\\
\end{array}
\right]
+
 \left[
\begin{array}{c}
x(k+1)\\
z(k)\\
\end{array}
\right]^T
 \left[
 \begin{array}{cc}
  P & 0\\
  0 & I_{\nz}\\
 \end{array}
 \right]
 \left[
\begin{array}{c}
x(k+1)\\
z(k)\\
\end{array}
\right]\le 0
\end{array}$},\\
(k=0,1,\cdots).  
\end{array}
\]
By summing up the above inequalities up to $k=N$, 
we have
\begin{equation}
\begin{array}{@{}l}
 x(N+1)^TPx(N+1)
-\tilgam^2 \sum_{k=0}^N |w(k)|_2^2\\
\hspace*{10mm}
+\sum_{k=0}^N w(k)^TQw(k)
+\sum_{k=0}^N |z(k)|_2^2\le 0.  
\end{array}
\label{eq:sum}
\end{equation}
We see $x(N+1)^TPx(N+1)\ge 0$ since $P\in\PSD_n$.  
On the other hand, since $w\in l_{2+}$ and since $Q\in \COP_{\nw}$, 
we see that $\sum_{k=0}^N w(k)^TQw(k)\ge 0$.  
Therefore, by letting $N\to\infty$ in \rec{eq:sum}, 
we have
$\|z\|_2^2 \le \tilgam^2 \|w\|_2^2=\tilgam^2$.  
Since this condition holds for arbitrary $w\in l_{2+}$ with $\|w\|_2=1$, 
we can conclude that
\[
 \begin{array}{@{}lcl}
 \|G\|_{2+}&=&\sup_{w\in l_{2+},\ \|w\|_2=1} \ \|z\|_2 \quad
  \le\quad \tilgam \quad < \quad \gamma.    
 \end{array}
\]
This completes the proof.  
\end{proofof}

On the basis of \rth{th:main}, let us consider the COP:
\begin{equation}
 \olgam_+:=\inf_{\gamma,P,Q}\ \mbox{subject to}\ \rec{eq:LMI},\ 
  P\in\PSD_n,\ 
  Q\in\COP_{\nw}.  
\label{eq:COP}
\end{equation}
In relation to this COP, recall that
\[
 \|G\|_2=\inf_{\gamma,P}\ \mbox{subject to}\ \rec{eq:LMI},\ 
  P\in\PSD_n,\ 
  Q=0.  
\]
It follows that $\|G\|_{2+}\le \olgam_+\le \|G\|_2$.  
Unfortunately, as we have already mentioned, 
it is hard to solve the COP \rec{eq:COP} in general.  
However, an upper bound of $\olgam_+$
can be computed efficiently by replacing $\COP$ in 
\rec{eq:COP} by $\PSD+\NN$ as follows:  
\begin{equation}
\begin{array}{@{}l}
 \dolgam_{+}:=\inf_{\gamma,P,Q}\ \mbox{subject to}\ \rec{eq:LMI},\\ 
  P\in\PSD_n,\ 
  Q\in\PSD_{\nw}+\NN_{\nw}.  
\end{array}
\label{eq:DNN}
\end{equation}
Note that this problem is essentially an SDP and hence tractable.  
We can readily see that 
$\|G\|_{2+}\le \olgam_+\le \dolgam_+\le \|G\|_2$ holds.    

Up to this point, we have described the basic idea of
the (upper bound) computation of $\|G\|_{2+}$.  
However, in the case where $\nw=1$, i.e., 
if the system $G$ has only a single disturbance input, 
then it is very clear that $\olgam_{+}=\|G\|_{2}$.  
This is because, since
$\COP_1=\PSD_1=\bbR_+$, and since the variable $Q$ enters
in block-diagonal part in \rec{eq:LMI}, 
we see that the optimal value of $Q$ in COP \rec{eq:COP} is zero.  
Namely, if $\nw=1$, 
it is impossible to obtain
an upper bound of $\|G\|_{2+}$ 
which is better than the trivial upper bound  
$\|G\|_{2}$ if we directly work on \rec{eq:COP}.  
In addition, 
we also deduce from this fact that 
the improvement of $\olgam_{+}$ over $\|G\|_{2}$
might not be significant if $G$ has a few number of 
disturbance inputs.  
To get around this difficulty, 
in the next section,    
we employ the discrete-time system lifting \cite{Bittanti_1996}.   

\subsection{Better Upper Bound Computation by System Lifting}
\label{sub:upper}

By applying the $N$-th order discrete-time lifting \cite{Bittanti_1996} 
to \rec{eq:G} with $N\in\bbN$, we can obtain another 
discrete-time LTI system $\hatG_N$ of the form
\begin{equation}
\hatG_N:\ 
\left\{
\arraycolsep=0.5mm
\begin{array}{cccccccc}
 \hatx(\kappa+1)&=& \hatA_N x(\kappa)& + &\hatB_N \hatw(\kappa), &\ \hatx(0)=0,\\
   \hatz(\kappa)&=& \hatC_N x(\kappa)& + &\hatD_N \hatw(\kappa)
\end{array}
\right.  
\label{eq:hatG}
\end{equation}
where
\begin{equation}
\scalebox{0.85}{$
\begin{array}{@{}l}
 \hatA_N :=A^N,\quad 
 \hatB_N:=
  \begin{bmatrix}
   A^{N-1}B & \cdots & AB & B 
  \end{bmatrix},\\
 \hatC_N:=
  \begin{bmatrix}
   C \\ CA \\ \vdots \\ CA^{N-1}
  \end{bmatrix},\quad 
  \hatD_N :=
  \begin{bmatrix}
   D          & 0      & \cdots & \cdots & 0 \\
   CB         & \ddots & \ddots & \ddots & \vdots \\
   CAB        & \ddots & \ddots & \ddots & \vdots \\
   \vdots     & \ddots & \ddots & \ddots & 0 \\
   CA^{N-1}B  & \cdots & CAB    & CB     & D
  \end{bmatrix}  
\end{array}$}
\label{eq:hABCD}
\end{equation}
and
\begin{equation}
\scalebox{0.9}{$
\begin{array}{@{}l}
 \hatw(\kappa)=
  \begin{bmatrix}
   w(\kappa N) \\
   w(\kappa N+1) \\
   \vdots\\
   w((\kappa+1)N-1) \\
  \end{bmatrix},\ 
 \hatz(\kappa)=
  \begin{bmatrix}
   z(\kappa N) \\
   z(\kappa N+1) \\
   \vdots\\
   z((\kappa+1)N-1) \\
  \end{bmatrix}.  
\end{array}$}
\label{eq:hzw}
\end{equation}
It is very clear that 
$G$ is stable (i.e., $A$ is Schur-Cohn stable) if and only if 
$\hatG_N$ is stable (i.e., $\hatA_N$ is Schur-Cohn stable).  
In addition, from \rec{eq:hzw}, 
we can readily see that 
$\|G\|_2=\|\hatG_N\|_2$ and
$\|G\|_{2+}=\|\hatG_N\|_{2+}$.  
With these facts in mind, for given $N\in\bbN$, let us define
\[
 \ol{\gamma}_{N+}:=\inf_{\gamma,P,Q}\ \mbox{subject to}\
\]
\begin{equation}
\begin{array}{@{}l}
\arraycolsep=0.5mm
L(\hatA_N,\hatB_N,\hatC_N,\hatD_N,P,Q,\gamma)\prec 0,
\end{array}
\label{eq:hLMI} 
\end{equation}
\[
  P\in\PSD_n,\ Q\in\COP_{N\nw},   
\]
\begin{equation}
\begin{array}{@{}l}
 \ol{\ol{\gamma}}_{N+}:=\inf_{\gamma,P,Q}\ \mbox{subject to}\
  \rec{eq:hLMI},\\ 
  P\in\PSD_n,\ 
  Q\in\PSD_{N\nw}+\NN_{N\nw}.  
\end{array}
\end{equation}
Then, we have $\forall N\in\bbN$ that 
\[
\|G\|_{2+}=\|\hatG_N\|_{2+}\le \ol{\gamma}_{N+}\le
\ol{\ol{\gamma}}_{N+}\le 
\|\hatG_N\|_{2}=\|G\|_2.  
\]
Namely, $\ol{\gamma}_{N+}$ and $\ol{\ol{\gamma}}_{N+}$ are upper bounds
of $\|G\|_{2+}$ and the latter is easy to compute.  
In particular, it is worth mentioning that 
we can obtain better (no worse) upper bounds by increasing $N$
as shown in the next theorem.  
The proof of this theorem is given at appendix section.  
\begin{theorem}
For given $N_1,N_2\in\bbN$ with $N_2=pN_1\ (\exists p\in\bbN)$, 
we have
\begin{equation}
\ol{\gamma}_{N_2+} \le \ol{\gamma}_{N_1+},\  
\ol{\ol{\gamma}}_{N_2+} \le \ol{\ol{\gamma}}_{N_1+}.   
\label{eq:up}
\end{equation}
\label{th:up}
\end{theorem}

We demonstrate the effectiveness of 
the lifting-based treatment in \rsub{sub:num}.  

\subsection{Lower Bound Computation by System Lifting}
\label{sub:lower}

In the preceding subsection, we consider the upper bound computation
of $\|G\|_{2+}$ for the discrete-time LTI system $G$.  
However, if we merely compute upper bounds, 
it is inherently impossible to evaluate their accuracy.  
To remedy this, in the section, we consider a method
for lower bound computation.   

In \rec{eq:hABCD}, recall that
the matrix $\hatD_N$ captures the input-output behavior
of the system $G$ up to time instant $N-1$.  Namely, we have
\[
  \begin{bmatrix}
   z(0) \\
   z(1) \\
   \vdots\\
   z(N-1) \\
  \end{bmatrix}=\hatD_N  
  \begin{bmatrix}
   w(0) \\
   w(1) \\
   \vdots\\
   w(N-1) \\
  \end{bmatrix}.  
\]
Therefore we can readily see that 
\begin{equation}
\|G\|_{2+}\ge \|\hatD_N\|_{2+}\ (\forall N\in\bbN).  
\label{eq:lb}
\end{equation}
It is also true that $\|\hatD_N\|_{2+}$ is monotonically 
non-decreasing with respect to $N\in\bbN$ and
$\|G\|_{2+}=\lim_{N\to\infty} \|\hatD_N\|_{2+}$.  
Therefore, if we can compute $\|\hatD_N\|_{2+}$ for each $N\in\bbN$ exactly, 
we can construct a monotonically nondecreasing
sequence of the lower bounds that converges to $\|G\|_{2+}$.  

As is well known, ``the maximal singular value'' $\|\hatD_N\|_{2}$
is characterized by the SDP:
\[
\|\hatD_N\|_{2}^2=\inf_{\gamsq}\ \gamsq\ \mbox{subject to}\ 
\gamsq I_{N\nw} - \hatD_N^T\hatD_N \succ 0.  
\]
Similarly, we see that
the ``positive'' maximal singular value $\|\hatD_N\|_{2+}$
is characterized by the COP:
\[
\|\hatD_N\|_{2+}^2=\inf_{\gamsq}\ \gamsq\ \mbox{subject to}\ 
\gamsq I_{N\nw} - \hatD_N^T\hatD_N \succ_{\COP} 0.  
\]
As repeatedly noticed, unfortunately, the above COP is 
numerically intractable if $N\nw>4$.  
A possible remedy is to replace the cone $\COP$ in the above COP 
by $\PSD+\NN$ and consider 
\begin{equation}
\inf_{\gamsq}\ \gamsq\ \mbox{subject to}\ \gamsq I_{N\nw} - \hatD_N^T\hatD_N
\succ_{\PSD+\NN} 0.    
\label{eq:PSD+NN}
\end{equation}
This is essentially an SDP.  
Nevertheless, it is quite important to note that the SDP 
\rec{eq:PSD+NN} provides an {\it upper bound} of 
(the square of)
$\|\hatD_N\|_{2+}$ and hence this is not fully fitted to  
our purpose here.  
Note that our goal here is to compute a lower bound
of $\|G\|_{2+}$ by way of \rec{eq:lb}
and hence what is required is to compute a 
{\it lower bound} of $\|\hatD_N\|_{2+}$.   

To compute a lower bound of $\|\hatD_N\|_{2+}$
in numerically tractable fashion, 
let us consider the dual of the SDP \rec{eq:PSD+NN} 
that is given as follows:
\begin{equation}
\begin{array}{@{}l}
\sup_{Z_N}\ \trace(\hatD_N^T\hatD_N Z_N)\ \mbox{subject to}\\ 
\trace(Z_N)=1,\quad Z_N\in\DNN_{N\nw}.  
\end{array}
\label{eq:PSDcapNN}
\end{equation}
Again, this is essentially an SDP and hence numerically tractable.  
Here, it is very clear that the (primal) SDP \rec{eq:PSD+NN}
has an interior point solution.  
Therefore, there is no duality gap between
the SDPs \rec{eq:PSD+NN} and \rec{eq:PSDcapNN}, and 
in particular the SDP \rec{eq:PSDcapNN} has an optimal solution
\cite{Klerk_2002}.  
It follows from the zero duality gap that the SDP \rec{eq:PSDcapNN} again 
provides an upper bound of $\|\hatD_N\|_{2+}$ and hence
we have to go further.  

With the above mentioned facts in mind, let us denote by 
$Z_N^\star$ an optimal solution of the SDP \rec{eq:PSDcapNN}.  
Moreover, let $v_N^\star\in\bbR_+^{N\nw}$ denote 
the unit eigenvector corresponding to the maximal eigenvalue 
of $Z_N^\star$.  
It should be noted that,  
from Perron-Frobenius theorem \cite{Horn_1985}, 
we can confirm that $v_N^\star$ is 
certainly nonnegative since $Z_N^\star\ge 0$.  
Then, if we define $\ul{\gamma}_N:=|\hatD_N v_N^\star|_2$, 
it is very clear that 
$\ul{\gamma}_N\le \|\hatD_N\|_{2+}$.  
This idea of lower bound computation comes from 
the rank-one exactness verification test
for LMI relaxation, which is 
frequently employed in the literature,  
see, e.g., \cite{Scherer_EJC2006}.  
Namely, it is straightforward to see that if 
$\rank(Z_N^\star)=1$ then $\ul{\gamma}_N= \|\hatD_N\|_{2+}$.  
We finally note that $\ul{\gamma}_N= \|\hatD_N\|_{2+}$
always holds if $\hatD_N^T\hatD_N\ge 0$.  
This is a direct consequence again from 
Perron-Frobenius theorem.  

\subsection{Exact Computation}
\label{sub:exact}

In some special cases, 
we can compute $\|G\|_{2+}$
exactly by solving an SDP.  
For instance, let us consider the case where
$G$ is ``static'' and its input-output property is given by
$z(k)=Dw(k)$.  
Then, we can readily see that 
\[
\|G\|_{2+}^2=\|D\|_{2+}^2=\inf_{\gamsq} \gamsq-D^TD\succeq_{\COP} 0.\\ 
\]
Since this COP is essentially 
an SDP if $\nw\le 4$, 
we arrive at the conclusion that we can compute 
$\|G\|_{2+}$ exactly by solving an SDP in the above special case.  

\subsection{Numerical Examples}
\label{sub:num}

Let us consider the case where
the coefficient matrices of the system
\rec{eq:G} are given by
\[
\begin{array}{@{}l}
 A=
\left[
 \begin{array}{rrrr}
 0.27 &  0.06 & -0.24 &  0.19 \\ 
-0.26 & -0.18 &  0.35 &  0.43 \\ 
 0.06 & -0.88 & -0.78 &  0.27 \\ 
-0.07 &  0.11 & -0.25 & -0.01 \\ 
 \end{array}\right],\ 
B = 
\left[
 \begin{array}{r}
  0.68 \\
  1.46 \\
 -0.22 \\
  0.45 \\
 \end{array}\right],\vspace*{2mm}\\ 
C = 
\left[
 \arraycolsep=2.2mm
 \begin{array}{rrrr}
  0.33 & -2.06 & 1.22 & 1.12
 \end{array}\right],\ 
D = 0.05.    
\end{array}
\]

\begin{figure}[t]
\vspace*{-3mm}
\begin{center}
\hspace*{-8mm}
 \includegraphics[scale=0.65]{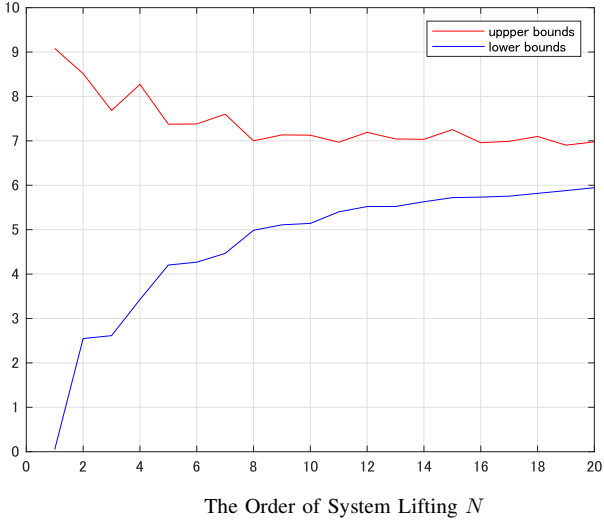}\vspace*{-4mm}\\
 \hspace*{8mm}\scalebox{0.8}{The Order of System Lifting $N$}
 \caption{Upper and Lower Bounds of $\|G\|_{2+}$.}
 \label{fig:lbub}
\end{center}
\vspace*{-5mm}
\end{figure}

By applying the discrete-time system lifting
and following the ideas in Subsections 
\ref{sub:upper} and \ref{sub:lower},    
we computed upper and lower bounds of $\|G\|_{2+}$.  
The results are shown in \rfig{fig:lbub}.  
The $l_2$ induced norm of $G$ turned out to be $\|G\|_2=9.0797$.  
On the other hand, the best upper bound obtained by lifting is
$\ol{\ol{\gamma}}_{19+}=6.9034$
and the best lower bound obtained by lifting is $\ul{\gamma}_{20}=5.9453$.  
The obtained upper bounds are NOT monotonically decreasing, 
but we can confirm that \rec{eq:up} is surely satisfied.  

\section{Stability Analysis of RNN with ReLU}
\label{sec:RNN}

In this section, we demonstrate the usefulness
of the ``positive'' $l_2$ induced norm 
in stability analysis of Recurrent Neural Networks (RNNs).  

\subsection{Basics of RNN and Stability}
\label{sub:RNN}

Let us consider the dynamics of 
the discrete-time RNNs typically described by
\begin{equation}
\left\{
\arraycolsep=0.5mm
 \begin{array}{ccl}
  x(k+1)&=&\Lambda x(k)+\Win w(k)+v(k),\\
  z(k)&=&\Wout x(k),\\
  w(k)&=&\Phi(z(k)+s(k))
 \end{array}
\right.  
\label{eq:RNN}
\end{equation}
where $x\in\bbR^n$ is the state and 
$\Lambda\in\bbR^{n\times n}$, 
$\Wout\in\bbR^{m\times n}$, 
$\Win\in\bbR^{n\times m}$ are constant matrices
with $\Lambda$ being Schur-Cohn stable. 
On the other hand, note that 
$s:\ [0,\infty)\to \bbR^m$ and $v:\ [0,\infty)\to \bbR^n$ are 
external input signals  
and $\Phi:\ \bbR^m\to \bbR^m$ is the
static activation function typically being nonlinear.  
The matrices $\Win$ and $\Wout$ are constructed from 
the weightings of the edges in RNN.  
We assume $x(0)=0$.  

In relation to \rec{eq:RNN}, let us define
\begin{equation}
\begin{array}{@{}l}
G:=
 \left[
 \begin{array}{c|cc}
  \Lambda & \Win & I_n\\ \hline
  \Wout & 0 & 0
 \end{array}
 \right],\\  
G_0:=
 \left[
 \begin{array}{c|c}
  \Lambda & \Win \\ \hline
  \Wout & 0 
 \end{array}
 \right],\quad  
G_1:=
 \left[
 \begin{array}{c|c}
  \Lambda & I_n\\ \hline
  \Wout & 0 
 \end{array}
 \right]  
\end{array}
\label{eq:Gs}
\end{equation}
where the initial states of these three systems are all zeros.  
Then the dynamics of the RNN given by 
\rec{eq:RNN} can be represented by the block-diagram 
shown in \rfig{fig:RNN}. 
We consider the typical case where the activation function is
Rectified Linear Unit (ReLU) whose input-output property is given by
\begin{equation}
\begin{array}{@{}l}
 \Phi(\xi)=\left[\ \phi(\xi_1)\ \cdots\ \phi(\xi_m)\ \right]^T,\\
 \phi: \bbR\to \bbR,\quad
 \phi(\eta)=
  \left\{
   \begin{array}{cc}
    \eta & (\eta\ge 0),  \\
    0 & (\eta< 0).  \\
   \end{array}
  \right.    
\end{array}  
\label{eq:ReLU}
\end{equation}
\begin{center}
\begin{figure}[b]
\begin{center}
\vspace*{-2mm}
\scalebox{1.0}{\begin{picture}(4.5,3)(0,0)
\put(0,2.5){\vector(1,0){0.85}}
\put(0.425,2.7){\makebox(0,0)[b]{$s$}}
\put(1,2.5){\circle{0.3}}
\put(0.7,2.7){\makebox(0,0){$ \scriptstyle + $}}
\put(1.15,2.5){\vector(1,0){0.85}}
\put(2,2){\framebox(1.5,1){$\Phi$}}
\put(3.5,2.5){\line(1,0){1.0}}
\put(4,2.7){\makebox(0,0)[b]{$w$}}
\put(4.5,2.5){\line(0,-1){1.75}}
\put(4.5,0.75){\vector(-1,0){1.0}}
\put(2,0){\framebox(1.5,1){$G$}}
\put(2,0.5){\line(-1,0){1.0}}
\put(1.5,0.3){\makebox(0,0)[t]{$z$}}
\put(1,0.5){\vector(0,1){1.85}}
\put(0.8,2.2){\makebox(0,0){$ \scriptstyle + $}}
\put(4.5,0.25){\vector(-1,0){1.0}}
\put(4,0.05){\makebox(0,0)[t]{$v$}}
\end{picture}}
\caption{Block-Diagram Representation of RNN.}
\label{fig:RNN}
\end{center}
\end{figure}
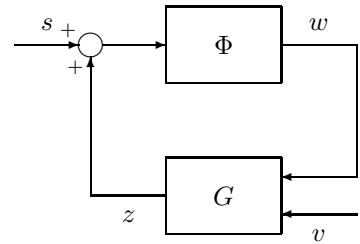
\end{center}

Since here we are dealing with nonlinear systems, 
it is of prime importance to clarify the definition of 
``stability.''
The definition we employ for the analysis of RNN is as follows.  
\begin{definition}\cite{Khalil_2002} (Finite Gain $l_2$ Stability) 
An operator $H:\  l_{2e}\ni u \to y\in l_{2e}$
is said to be finite gain $l_2$ stable
if there exists a non-negative constant $\gamma$ such that
$\|y_\tau\|_2 \le \gamma \|u_\tau\|_2$ holds
for any $u \in l_{2e}$ and $\tau \in [0,\infty)$. 
\end{definition}

In the following, we analyze the 
finite gain $l_2$ stability of the operator in RNN
shown in \rfig{fig:RNN} with respect to input 
$[\ s^T\ v^T\ ]^T\in l_{2e}$
and output $[\ z^T\ w^T\ ]^T \in l_{2e}$.  
Note that the feedback connection in \rfig{fig:RNN}
is well-posed since its dynamics is given by the state-space equation
\rec{eq:RNN}.  
We also note that we implicitly use the causality of $G$ and $\Phi$.

\subsection{Small Gain Type Stability Conditions}

The next theorem provides a 
small gain type stability condition with 
``positive'' $l_2$ induced norm.  
\begin{theorem}
The RNN given by \rec{eq:RNN} with the activation function ReLU given by \rec{eq:ReLU}
is finite-gain $l_2$ stable
if $\|G_0\|_{2+}< 1$ where $G_0$ is given by \rec{eq:Gs}.  
\label{th:SG_RNN}
\end{theorem}
\begin{proofof}{\rth{th:SG_RNN}}
We can readily see that $\|\Phi\|_2=1$.  
It is also true from \rec{eq:ReLU} that 
$w=\Phi(z+s)$ is nonnegative whatever 
$z$ and $s$ are.  
Therefore, for any $[\ s^T\ v^T\ ]^T\in l_{2e}$ 
and $\tau\in[0,\infty)$, we have
\begin{equation}
\begin{array}{@{}lcl}
 \|z_\tau\|_{2} &= & 
  \left\|\left(G\left[
	  \begin{array}{c}
	   w\\
	   v
	  \end{array}\right]\right)_\tau\right\|_2\vspace*{2mm}\\ 
 &\le & 
  \left\|\left(G\left[
	  \begin{array}{c}
	   w\\
	   0
	  \end{array}\right]\right)_\tau\right\|_2+
  \left\|\left(G\left[
	  \begin{array}{c}
	   0\\
	   v
	  \end{array}\right]\right)_\tau\right\|_2\\
 &=&  \|\left(G_0\Phi (z+s)\right)_\tau\|_2+\|(G_1v)_\tau\|_2\\
 &\le&  \|G_0\Phi (z+s)_\tau\|_2+\|G_1v_\tau\|_2\\
 &=&  \|G_0\Phi (z_\tau+s_\tau)\|_2+\|G_1v_\tau\|_2\\
 &\le&  \|G_0\|_{2+}\|\Phi\|_2(\|z_\tau\|_2+\|s_\tau\|_2)+\|G_1\|_2\|v_\tau\|_2\hspace*{-10mm}\\
 &=&  \|G_0\|_{2+}\|z_\tau\|_2+\|G_0\|_{2+}\|s_\tau\|_2+\|G_1\|_2\|v_\tau\|_2,\hspace*{-10mm}\\
\end{array}
\label{eq:ztau}
\end{equation}
\begin{equation}
\begin{array}{@{}lcl}
 \|w_\tau\|_{2} &=& \|(\Phi (z+s))_\tau\|_2\\
 &\le&  \|\Phi\|_2 (\|z_\tau\|_2+\|s_\tau\|_2)\\
 &=&  \|z_\tau\|_2+\|s_\tau\|_2.  
\end{array}
\label{eq:wtau}
\end{equation}
If we define 
\[
\gamma_0:=\|G_0\|_{2},\  
\gamma_{0+}:=\|G_0\|_{2+},\ 
\gamma_{1}:=\|G_1\|_{2}
\]
and assume $\gamma_{0+}<1$, then we readily obtain
\[
\begin{array}{@{}lcl}
\|z_\tau\|_{2} \le (1-\gamma_{0+})^{-1} \gamma_{0+} \|s_\tau\|_2
+(1-\gamma_{0+})^{-1} \gamma_1 \|v_\tau\|_2,\\ 
\|w_\tau\|_{2} \le (1-\gamma_{0+})^{-1} \|s_\tau\|_2
+(1-\gamma_{0+})^{-1} \gamma_1 \|v_\tau\|_2.  
\end{array}
\]
It follows from \rle{le:ineq} given in the appendix section that
\[
\scalebox{0.86}{$
\begin{array}{@{}lcl}
\left\|
\left[
\begin{array}{c}
 z_\tau\\
 w_\tau
\end{array}\right]
\right\|_2\le \sqrt{2}
\left\|
\left[
\begin{array}{cc}
(1-\gamma_{0+})^{-1} \gamma_{0+} & (1-\gamma_{0+})^{-1} \gamma_1\\ 
(1-\gamma_{0+})^{-1} & (1-\gamma_{0+})^{-1} \gamma_1
\end{array}\right]
\right\|_2
\left\|
\left[
\begin{array}{c}
 s_\tau\\
 v_\tau
\end{array}\right]
\right\|_2  
\end{array}$}
\]
holds for all $[\ s^T\ v^T\ ]^T\in l_{2e}$ and $\tau\in[0,\infty)$.  
Therefore we can conclude that
RNN given by \rec{eq:RNN} with ReLU given by \rec{eq:ReLU}
is finite-gain $l_2$ stable
if $\gamma_{0+}=\|G_0\|_{2+}< 1$.  
\end{proofof}

The stability condition $\|G_0\|_{2+}<1$ in \rth{th:SG_RNN}
is of course 
a milder condition than the
``standard'' small gain condition that requires $\|G_0\|_{2}< 1$.  
From \rth{th:main}, 
we see that 
$\|G_0\|_{2+}< 1$ holds if 
there exist $P\in\PSD_n$ and $Q\in\COP_{m}$
such that
\[
\scalebox{0.95}{$
\begin{array}{@{}l}
\arraycolsep=0.5mm
\left[
\begin{array}{cc}
 -P & 0 \\
 0 & -I_{m}+Q
\end{array}
\right]+
 \left[
 \begin{array}{cc}
  \Lambda & \Win \\
  \Wout & 0\\
 \end{array}
 \right]^T
 \left[
 \begin{array}{cc}
  P & 0\\
  0 & I_{m}\\
 \end{array}
 \right]
 \left[
 \begin{array}{cc}
  \Lambda & \Win \\
  \Wout & 0\\
 \end{array}
 \right]\prec 0.  
\end{array}$}
\]
%

\subsection{Scaled Small Gain Type Stability Condition}

It is not hard to see that ReLU $\Phi$ satisfies
$\Phi(\xi)=(D^{-1}\Phi D)(\xi)$ for any 
$D\in\bbD_{++}^{n}$ where
$\bbD_{++}^n\subset\bbR^{n\times n}$ stands for the set of diagonal matrices
with strictly positive diagonal entries.  
Therefore we readily deduce that the RNN
given by \rec{eq:RNN} with ReLU given by \rec{eq:ReLU}
is finite-gain $l_2$ stable
if there exists $D\in\bbD_{++}^n$ such that
$\|D^{-1}G_0D\|_{2+}< 1$.  
From \rth{th:main}, this condition holds if
there exist $P\in\PSD_n$, $D\in\bbD_{++}^n$, and $Q\in\COP_{m}$
such that
\[
\scalebox{0.82}{$
\begin{array}{@{}l}
\arraycolsep=0.1mm
\left[
\begin{array}{cc}
 -P & 0 \\
 0 & -I_m+Q
\end{array}
\right]+
 \left[
 \begin{array}{cc}
  \Lambda & \Win D \\
  D^{-1}\Wout & 0\\
 \end{array}
 \right]^T
 \left[
 \begin{array}{cc}
  P & 0\\
  0 & I_n\\
 \end{array}
 \right]
 \left[
 \begin{array}{cc}
  \Lambda & \Win D \\
  D^{-1}\Wout & 0\\
 \end{array}
 \right]\prec 0.  
\end{array}$}
\]
We can equivalently translate this nonconvex condition
to the convex condition that
there exist $P\in\PSD_n$, $S\in\bbD_{++}^m$, and $Q\in\COP_m$
such that
\begin{equation}
\scalebox{0.90}{$
\begin{array}{@{}l}
\arraycolsep=0.5mm
\left[
\begin{array}{cc}
 -P & 0 \\
 0 & -S+Q
\end{array}
\right]+
 \left[
 \begin{array}{cc}
  \Lambda & \Win \\
  \Wout & 0\\
 \end{array}
 \right]^T
 \left[
 \begin{array}{cc}
  P & 0\\
  0 & S\\
 \end{array}
 \right]
 \left[
 \begin{array}{cc}
  \Lambda & \Win\\
  \Wout & 0\\
 \end{array}
 \right]\prec 0.  
\end{array}$}
\label{eq:RNNSGSLMI}
\end{equation}
We note that the 
``standard'' small gain condition is recovered 
if we let $Q=0$ in \rec{eq:RNNSGSLMI}.  


\subsection{Numerical Examples}

In \rec{eq:RNN}, let us consider the case 
$\Lambda=0$, $\Wout=I_6$ and 
\[
\scalebox{1.0}{$
\begin{array}{@{}l}
\Win= 
\left[
\begin{array}{rrrrrr}
 0.29 & -0.04 &  0.02+a & -0.35 & -0.05 & -0.12\\
-0.29 & -0.24 & -0.01 &  0.12 & -0.13 &  0.18\\
-0.50 &  b &  0.23 &  0.40 & -0.28 & -0.08\\
 0.14 & -0.27 & -0.15 &  0.13 & -0.47 & -0.28\\
-0.10 & -0.10 &  0.08 &  0.14 & -0.22 &  0.50\\
-0.11 & -0.28 & -0.21 & -0.14 & -0.09 &  0.20\\
\end{array}
\right].  
\end{array}$}
\]
For $(a,b)=(0,0)$ we see $\|G_0\|_2=0.9605$.  
Here we examined the finite-gain $l_2$ stability 
over the (time-invariant) parameter variation 
$a\in[-8,8]$ and $b\in [-8,8]$.  
We tested the following stability conditions.    

\noindent
\begin{equation}
\begin{array}{@{}l}
\mbox{
SSG:}\\
\mbox{
Find $P\in \PSD_n$, $S\in\bbD_{++}^m$}\\
\mbox{
such that \rec{eq:RNNSGSLMI} holds where $Q=0$}.  
\end{array}
\label{eq:SSG}
\end{equation}
\begin{equation}
\begin{array}{@{}l}
\mbox{
SSG+COP:}\\
\mbox{
Find $P\in \PSD_n$, $S\in\bbD_{++}^m$, 
$Q\in\PSD_m+\NN_m$}\\
\mbox{
such that \rec{eq:RNNSGSLMI} holds}.  
\end{array}
\label{eq:SSG+COP}
\end{equation}

\begin{figure}[t]  
\begin{center}
 \vspace*{-3mm}
 \hspace*{-5mm}
 \includegraphics[scale=0.65]{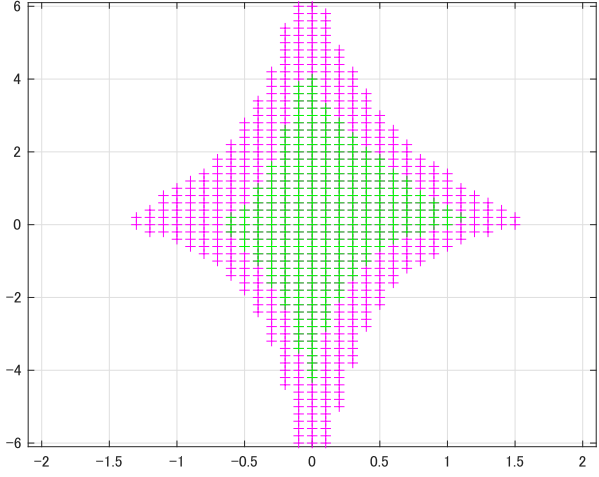}
 \vspace*{-12mm}
 \caption{Comparison: SSG vs SSG+COP.  }
 \label{fig:comp1}
 \vspace*{-8mm}
\end{center}
\end{figure}

It is very clear that 
if \rec{eq:SSG} is feasible then \rec{eq:SSG+COP} is.  
In \rfig{fig:comp1}, 
we plot $(a,b)$ for which the RNN is proved to be stable
by the above stability conditions.  
Both LMIs \rec{eq:SSG} and \rec{eq:SSG+COP} turned out to be feasible 
for $(a,b)$ with green plot, 
whereas only \rec{eq:SSG+COP} turned out to be feasible
for $(a,b)$ with magenta plot.   
We can clearly see the effectiveness
of the present new stability condition with
the ``positive'' $l_2$ induced norm.  

\section{Conclusion and Future Works}

In this paper, 
we newly introduced the ``positive'' $l_2$ induced norm of 
discrete-time LTI systems
where the input signals are restricted to be nonnegative.  
On the basis of copositive programming, 
we provide tractable methods for the upper and lower bound computation
of the ``positive'' $l_2$ induced norm.   
We illustrate its usefulness in stability analysis 
of recurrent neural networks with activation functions
being rectified linear units.  

The present paper just described
basic treatments of the ``positive'' $l_2$ induced norm
and its application.  
In closing, 
we summarize outstanding issues to be investigated.  

\subsubsection{Treatment of COP}
In the present paper, we converted a COP to an SDP
by simply replacing $\COP$ by $\PSD+\NN$.  
However, this treatment is primitive and hence
conservative.  In this respect, 
Lasserre \cite{Lasserre_2014} has already shown how to
construct a hierarchy of SDPs to solve COP 
in an asymptotically exact fashion.  
Nevertheless, this approach does not allow us to 
handle practical size problems 
since the size of SDPs grows very rapidly.  
We need further effort to reduce computational burden
for instance by finding out sparsity structure.  
We plan to rely on efficient first-order methods to 
solve the specific conic relaxations arising from 
polynomial optimization problems with sphere constraints
\cite{Mai_arXiv2020}.  

\subsubsection{Stability Analysis of Lurye Systems with COP Multipliers}

Our (scaled) small-gain type treatment for the stability analysis of
RNN might be too shallow in view of 
advanced integral quadratic constraint (IQC) theory \cite{Megretski_IEEE1997}.  
Namely, for the stability analysis of feedback systems 
constructed from an LTI system and 
nonlinear elements
(i.e., Lurye systems), the effectiveness of 
the IQC approach with Zames-Falb multipliers \cite{Zames_SIAM1968}
is widely recognized, see, e.g., \cite{Fetzler_IFAC2017,Fetzler_IJRN2017}.  
Therefore it is strongly preferable if we can build a new COP-based approach
on the basis of powerful IQC-based framework.  
To this end, we need to explore sound ways 
to capture the properties of
nonlinear elements exhibiting positivity (such as ReLU) 
by introducing copositive multipliers and incorporate
them into existing IQC conditions.   
It is also important to seek for possible ways 
to introduce copositive multipliers to deal with 
saturated systems 
on the basis of the techniques developed for 
their analysis and synthesis \cite{Tarbouriech_2011}.  

These topics are currently under investigation.  

%


\begin{thebibliography}{10}

\bibitem{Khalil_2002}
H.~Khalil.
\newblock {\em Nonlinear Systems}.
\newblock Prentice Hall, 2002.

\bibitem{Rantzer_SCL1996}
A.~Rantzer.
\newblock On the {K}alman-{Y}akubovich-{P}opov lemma.
\newblock {\em Systems and Control Letters}, 28(1):7--10, 1996.

\bibitem{Ebihara_IEEE2017}
Y.~Ebihara, D.~Peaucelle, and D.~Arzelier.
\newblock Analysis and synthesis of interconnected positive systems.
\newblock {\em IEEE Transactions on Automatic Control}, 62(2):652--667, 2017.

\bibitem{Kato_LCSS2020}
T.~Kato, Y.~Ebihara, and T.~Hagiwara.
\newblock Analysis of positive systems using copositive programming.
\newblock {\em IEEE Control Systems Letters}, 4(2):444--449, 2020.

\bibitem{Duer_2010}
M.~Duer.
\newblock Copositive programming - a survey.
\newblock In M.~Diehl, F.~Glineur, E.~Jarlebring, and W.~Michiels, editors,
  {\em Recent Advances in Optimization and Its Applications in Engineering},
  pages 3--20. Springer, 2010.

\bibitem{Barabanov_IEEENN2002}
N.~E. Barabanov and D.~V. Prokhorov.
\newblock Stability analysis of discrete-time recurrent neural networks.
\newblock {\em IEEE Transactions on Neural Networks}, 13(2):292--303, 2002.

\bibitem{Zhang_IEEENN2014}
H~Zhang, Z~Wang, and D.~Liu.
\newblock A comprehensive review of stability analysis of continuous-time
  recurrent neural networks.
\newblock {\em IEEE Transactions on Neural Networks and Learning Systems},
  25(7):1229--1262, 2014.

\bibitem{Salehinejad_2018}
H.~Salehinejad, S.~Sankar, J.~Barfett, E.~Colak, and S.~Valaee.
\newblock Recent advances in recurrent neural networks.
\newblock arXiv:1801.01078v3 [cs.NE], 2018.

\bibitem{Boyd_2004}
S.~Boyd and L.~Vandenberghe.
\newblock {\em Convex Optimization}.
\newblock Cambridge University Press, Cambridge, 2004.

\bibitem{Berman_1979}
A.~Berman and R.~J. Plemmons.
\newblock {\em Nonnegative Matrices in the Mathematical Sciences}.
\newblock Academic Press, New York, 1979.

\bibitem{Berman_2003}
A.~Berman and N.~Shaked-Monderer.
\newblock {\em Completely positive matrices}.
\newblock World Scientific Publishing, Singapore, 2003.

\bibitem{Duer_EJLA2008}
M.~Duer and G.~Still.
\newblock Interior points of the completely positive cone.
\newblock {\em Electronic Journal of Linear Algebra}, 17:48--53, 2008.

\bibitem{Dickinson_EJLA2010}
P.~J.~C. Dickinson.
\newblock An improved characterization of the interior of the completely
  positive cone.
\newblock {\em Electronic Journal of Linear Algebra}, 20:723--729, 2010.

\bibitem{Farina_2000}
L.~Farina and S.~Rinaldi.
\newblock {\em Positive Linear Systems: Theory and Applications}.
\newblock John Wiley and Sons, Inc., 2000.

\bibitem{Tanaka_IEEE2011}
T.~Tanaka and C.~Langbort.
\newblock The bounded real lemma for internally positive systems and
  {$H_\infty$} structured static state feedback.
\newblock {\em IEEE Transactions on Automatic Control}, 56(9):2218--2223, 2011.

\bibitem{Rantzer_IEEE2016}
A.~Rantzer.
\newblock On the {K}alman-{Y}akubovich-{P}opov lemma for positive systems.
\newblock {\em IEEE Transactions on Automatic Control}, 61(5):1346--1349, 2016.

\bibitem{Bittanti_1996}
S.~Bittanti and P.~Colaneri.
\newblock Analysis of discrete-time linear periodic systems.
\newblock In Cornelius~T. Leondes, editor, {\em Control and Dynamic Systems},
  volume~78, pages 313--339. Academic Press, New York, 1996.

\bibitem{Klerk_2002}
E.~Klerk.
\newblock {\em Aspects of Semidefinite Programming}.
\newblock Kluwer Academic Publishers, 2002.

\bibitem{Horn_1985}
R.~A. Horn and C.~A. Johnson.
\newblock {\em Matrix Analysis}.
\newblock Cambridge University Press, New York, 1985.

\bibitem{Scherer_EJC2006}
C.~W. Scherer.
\newblock {LMI} relaxations in robust control.
\newblock {\em European Journal of Control}, 12(1):3--29, 2006.

\bibitem{Lasserre_2014}
J.~B. Lasserre.
\newblock New approximations for the cone of copositive matrices and its dual.
\newblock {\em Mathematical Programming, Series A}, 144:265--276, 2014.

\bibitem{Mai_arXiv2020}
N.~H.~A. Mai, V.~Magron, and J.~B. Lasserre.
\newblock A hierarchy of spectral relaxations for polynomial optimization.
\newblock In {\em arXiv:2007.09027v1 [math.OC]}, 2020.

\bibitem{Megretski_IEEE1997}
A.~Megretski and A.~Rantzer.
\newblock System analysis via integral quadratic constraints.
\newblock {\em IEEE Transactions on Automatic Control}, 42(6):819--830, 1997.

\bibitem{Zames_SIAM1968}
G.~Zames and P.~Falb.
\newblock Stability conditions for systems with monotone and slope-restricted
  nonlinearities.
\newblock {\em SIAM Journal on Control}, 6(1):89--108, 1968.

\bibitem{Fetzler_IFAC2017}
M.~Fetzer and C.~W. Scherer.
\newblock Absolute stability analysis of discrete time feedback
  interconnections.
\newblock {\em IFAC PapersOnline}, 50(1):8447--8453, 2017.

\bibitem{Fetzler_IJRN2017}
M.~Fetzer and C.~W. Scherer.
\newblock Full-block multipliers for repeated, slope-restricted scalar
  nonlinearities.
\newblock {\em International Journal of Robust and Nonlinear Control},
  27(17):3376--3411, 2017.

\bibitem{Tarbouriech_2011}
S.~Tarbouriech, G.~Garcia, G.~Silva Jr., and I.~Queinnec.
\newblock {\em Stability and Stabilization of Linear Systems with Saturating
  Actuators}.
\newblock Springer, 2011.

\bibitem{Oliveira_IJC2002}
M.~C. de~Oliveira, J.~C. Geromel, and J.~Bernussou.
\newblock Extended {$H_2$} and {$H_\infty$} norm characterizations and
  controller parametrizations for discrete-time systems.
\newblock {\em International Journal of Control}, 75:666--679, 2002.

\end{thebibliography}

\appendix
\section*{Proof of \rth{th:up}}

For the proof we need the next two lemmas.  
\begin{lemma}
For the system $G$ given by \rec{eq:G} and 
given $N_1,N_2\in\bbN$, let us define
$(\hatA_{N_1},\hatB_{N_1},\hatC_{N_1},\hatD_{N_1})$ and
$(\hatA_{N_2},\hatB_{N_2},\hatC_{N_2},\hatD_{N_2})$ by
\rec{eq:hABCD}.  Then, we have
\[
 \begin{array}{@{}l}
 \hatA_{N_2}\hatA_{N_1}=\hatA_{N_1+N_2},\ 
  \left[\ \hatA_{N_2}\hatB_{N_1}\ \hatB_{N_2}\ \right]=\hatB_{N_1+N_2},\\
  \left[
   \begin{array}{c}
    \hatC_{N_1} \\
    \hatC_{N_2}\hatA_{N_1} \\
   \end{array}
  \right]=\hatC_{N_1+N_2},\ 
  \left[
   \begin{array}{cc}
    \hatD_{N_1} & 0 \\
    \hatC_{N_2}\hatB_{N_1} & \hatD_{N_2} 
   \end{array}
  \right]=\hatD_{N_1+N_2}.  
 \end{array}
\]
\label{le:1}
\end{lemma}
\begin{lemma}
For given 
$A_1\in\bbR^{n\times n}$, 
$B_1\in\bbR^{n\times m_1}$, 
$C_1\in\bbR^{l_1\times n}$, 
$D_1\in\bbR^{l_1\times m_1}$ and 
$A_2\in\bbR^{n\times n}$, 
$B_2\in\bbR^{n\times m_2}$, 
$C_2\in\bbR^{l_2\times n}$, 
$D_2\in\bbR^{l_2\times m_2}$ and 
$\gamma>0$, 
suppose there exist
$P\in\bbS_n$, $Q_1\in\bbS_{m_1}$ and 
$Q_2\in\bbS_{m_2}$ such that
\begin{equation}
 L(A_1,B_1,C_1,D_1,P,Q_1,\gamma)\prec 0,\\
\label{eq:LMIc1}
\end{equation}
\begin{equation}
 L(A_2,B_2,C_2,D_2,P,Q_2,\gamma)\prec 0.  
\label{eq:LMIc2}
\end{equation}
Then we have
\begin{equation}
\begin{array}{@{}l}
 L(\clA,\clB,\clC,\clD,P,\clQ,\gamma)\prec 0, 
\end{array}
\label{eq:LMIc12}
\end{equation}
\begin{equation}
\arraycolsep=0.5mm
 \begin{array}{@{}l}
 \clA:=A_2A_1,\ 
  \clB:=
  \left[\ A_{2}B_{1}\ B_{2}\ \right],\\
  \clC:=
  \left[
   \begin{array}{c}
    C_{1} \\
    C_{2}A_{1} \\
   \end{array}
  \right],\ 
 \clD:=
  \left[
   \begin{array}{cc}
    D_{1} & 0 \\
    C_{2}B_{1} & D_{2} 
   \end{array}
  \right],\
  \clQ:=
   \left[
    \begin{array}{cc}
     Q_1 & 0 \\
     0   & Q_2
    \end{array}
   \right].  
 \end{array}
\label{eq:clM}
\end{equation}
\label{le:2}
\end{lemma}

We can confirm the validity of \rle{le:1} by direct calculation.  
The proof of \rle{le:2} is given as follows.  

\begin{proofof}{\rle{le:2}}
From \cite{Oliveira_IJC2002}, 
we see that \rec{eq:LMIc1} holds if and only if
there exists $G_1\in \bbR^{n\times n}$ such that
\[
 L_e(A_1,B_1,C_1,D_1,P,Q_1,G_1,\gamma)\prec 0
\]
\[
\scalebox{0.83}{$
\begin{array}{@{}l}
L_e(A_1,B_1,C_1,D_1,P,Q_1,G_1,\gamma):=\\ 
\arraycolsep=0.5mm
\left[
\begin{array}{ccc}
 -P+C_1^TC_1 & C_1^TD_1 & 0 \\
 D_1^TC_1 & D_1^TD_1-\gamma^2 I_{m_1}+Q_1 &  0 \\
 0 & 0 & P
\end{array}
\right]+
\He\left\{
 \left[
 \begin{array}{c}
  A_1^T \\
  B_1^T \\
  -I
 \end{array}
 \right]
 \left[
 \begin{array}{ccc}
  0 & 0 & G_1
 \end{array}
 \right]\right\}.  
\end{array}$}
\]
Similarly, 
\rec{eq:LMIc2} holds if and only if
there exists $G_2\in \bbR^{n\times n}$ such that
$
 L_e(A_2,B_2,C_2,D_2,P,Q_2,G_2,\gamma)\prec 0.  
$
It follows from \rec{eq:LMIc1} and \rec{eq:LMIc2} that
\[
\scalebox{0.75}{$
\begin{array}{@{}l}
\left[
\begin{array}{cc}
 L_e(A_1,B_1,C_1,D_1,P,Q_1,G_1,\gamma)& 0 \\
 0 & L_e(A_2,B_2,C_2,D_2,P,Q_2,G_2,\gamma)\\
\end{array}
\right]\prec 0.  
\end{array}$}
\]
By multiplying the above inequality by
\[
\scalebox{1.0}{$
\begin{array}{@{}l}
\left[
 \begin{array}{cccccc}
  I_n & 0 & 0 & 0 & 0 & 0 \\
  0 & I_{m_1} & 0 & 0 & 0 & 0 \\
  0 & 0 & I_n & I_n & 0 & 0 \\
  0 & 0 & 0 & 0 & I_{m_2} & 0 \\
  0 & 0 & 0 & 0 & 0 & I_n \\
 \end{array}\right]
\end{array}$}
\]
from left and its transpose from right, we have
\begin{equation}
\scalebox{0.75}{$
\begin{array}{@{}l}
\arraycolsep=0.5mm
\left[
\begin{array}{ccccc}
 -P+C_1^TC_1 & C_1^TD_1 & 0 & 0 & 0 \\
 \ast & D_1^TD_1-\gamma^2 I_{m_1}+Q_1 & 0 & 0 & 0 \\
 \ast & \ast & C_2^TC_2 & C_2^TD_2 & 0 \\
 \ast & \ast & \ast & D_2^TD_2-\gamma^2 I_{m_2}+Q_2 & 0 \\
 \ast & \ast & \ast & \ast & P
\end{array}
\right]\\
+\He\left\{
 \left[
 \begin{array}{cc}
  A_1^T & 0 \\
  B_1^T & 0 \\
  -I    & A_2^T \\
  0     & B_2^T \\
  0     & -I
 \end{array}
 \right]
 \left[
 \begin{array}{ccccc}
  0 & 0 & G_1 & 0 & 0 \\
  0 & 0 & 0   & 0 & G_2
 \end{array}
 \right]\right\}\prec 0.  
\end{array}$}
\label{eq:connect}
\end{equation}
Since 
\[
\scalebox{0.9}{$
\begin{array}{@{}l}
  \left[
 \begin{array}{cc}
  A_1^T & 0 \\
  B_1^T & 0 \\
  -I    & A_2^T \\
  0     & B_2^T \\
  0     & -I
 \end{array}
 \right]^\perp
=
\left[
\begin{array}{ccccc}
 I_n & 0 & A_1^T & 0 & A_1^TA_2^T\\
 0 & I_{m_1} & B_1^T & 0 & B_1^TA_2^T\\
 0 & 0  & 0 & I_{m_2} & B_2^T
\end{array}
\right]=:J, 
\end{array}$}
\]
the inequality \rec{eq:connect} implies
\[
\scalebox{0.7}{$
\begin{array}{@{}l}
\arraycolsep=0.5mm
J \left[
\begin{array}{ccccc}
 -P+C_1^TC_1 & C_1^TD_1 & 0 & 0 & 0 \\
 \ast & D_1^TD_1-\gamma^2 I_{m_1}+Q_1 & 0 & 0 & 0 \\
 \ast & \ast & C_2^TC_2 & C_2^TD_2 & 0 \\
 \ast & \ast & \ast & D_2^TD_2-\gamma^2 I_{m_2}+Q_2 & 0 \\
 \ast & \ast & \ast & \ast & P
\end{array}
\right]J^T \prec 0
\end{array}$}
\]
or equivalently, 
\[
\scalebox{0.8}{$
\begin{array}{@{}l}
\left[
\begin{array}{ccc}
 -P  & 0 & 0 \\
 \ast & -\gamma^2 I_{m_1}+Q_1 & \\
 0 & 0 & -\gamma^2 I_{m_2}+Q_2  \\
\end{array}
\right]+\left[
  \begin{array}{c}
   C_1^T \\
   D_1^T \\
   0
  \end{array}
 \right]
\left[
  \begin{array}{c}
   C_1^T \\
   D_1^T \\
   0
  \end{array}
 \right]^T\\
+
\left[
  \begin{array}{c}
   A_1^TC_2^T \\
   B_1^TC_2^T \\
   D_2^T \\
  \end{array}
 \right]
\left[
  \begin{array}{c}
   A_1^TC_2^T \\
   B_1^TC_2^T \\
   D_2^T \\
  \end{array}
 \right]^T+
\left[
\begin{array}{c}
A_1^TA_2^T\\
B_1^TA_2^T\\
B_2^T
\end{array}
\right]P
\left[
\begin{array}{c}
A_1^TA_2^T\\
B_1^TA_2^T\\
B_2^T
\end{array}
\right]^T \prec 0.    
\end{array}$}
\]
From \rec{eq:clM}, this can be rewritten equivalently as
\[
\scalebox{0.85}{$
\arraycolsep=0.5mm
\begin{array}{@{}l}
\left[
\begin{array}{cc}
 -P  & 0 \\
 \ast & -\gamma^2 I_{m_1+m_2}+\clQ \\
\end{array}
\right]+\left[
  \begin{array}{c}
   \clC^T \\
   \clD^T \\
  \end{array}
 \right]
\left[
  \begin{array}{c}
   \clC^T \\
   \clD^T \\
  \end{array}
 \right]^T
+
\left[
\begin{array}{c}
\clA^T\\
\clB^T
\end{array}
\right]P
\left[
\begin{array}{c}
\clA^T\\
\clB^T
\end{array}
\right]^T \prec 0.    
\end{array}$}
\]
This clearly shows that \rec{eq:LMIc12} holds.  
\end{proofof}

We are now ready to prove \rth{th:up}.  

\begin{proofof}{\rth{th:up}}
We prove $\ol{\gamma}_{N_2+} \le \ol{\gamma}_{N_1+}$.  
The proof of $\ol{\ol{\gamma}}_{N_2+} \le \ol{\ol{\gamma}}_{N_1+}$
follows similarly.  
For the proof of 
$\ol{\gamma}_{N_2+} \le \ol{\gamma}_{N_1+}$,   
it suffices to show that if there exist
$P_{N_1}\in\PSD_n$ and $Q_{N_1}\in\COP_{N_1\nw}$ such that
\begin{equation}
L(\hatA_{N_1},\hatB_{N_1},\hatC_{N_1},\hatD_{N_1},P_{N_1},Q_{N_1},\gamma)
\prec 0 
\label{eq:N1} 
\end{equation}
for given $\gamma >0$, then 
there exist
$P_{N_2}\in\PSD_n$ and $Q_{N_2}\in\COP_{N_2\nw}$ such that
\begin{equation}
L(\hatA_{N_2},\hatB_{N_2},\hatC_{N_2},\hatD_{N_2},P_{N_2},Q_{N_2},\gamma)
\prec 0.  
\label{eq:N2} 
\end{equation}
To this end, we first note from \rec{eq:N1} and \rle{le:2} that
\[
\scalebox{0.63}{$
\begin{array}{@{}l}
L\left(\hatA_{N_1}^2,[\ \hatA_{N_1}\hatB_{N_1}\ \hatB_{N_1}\ ],
\left[
 \begin{array}{c}
  \hatC_{N_1}\\
  \hatC_{N_1}\hatA_{N_1}\\
 \end{array}
 \right], 
\left[
 \begin{array}{cc}
  \hatD_{N_1} & 0 \\
  \hatC_{N_1}\hatB_{N_1}& \hatD_{N_1}\\
 \end{array}
 \right], 
P_{N_1},
\left[
 \begin{array}{cc}
  Q_{N_1} & 0 \\
  0 & Q_{N_1}
 \end{array}
 \right], \gamma\right) 
\prec 0
\end{array}$}
\]
holds.  From \rle{le:1}, this can be rewritten equivalently as
\begin{equation}
\begin{array}{@{}l}
L\left(\hatA_{2N_1},\hatB_{2N_1},\hatC_{2N_1},\hatD_{2N_1},
P_{2N_1},Q_{2N_1}, \gamma\right) 
\prec 0,\\
P_{2N_1}:=P_{N_1},\ 
Q_{2N_1}:=\left[
 \begin{array}{cc}
  Q_{N_1} & 0 \\
  0 & Q_{N_1}
 \end{array}
 \right].  
\end{array}
\label{eq:2N1}
\end{equation}
Similarly, from \rec{eq:N1} and \rec{eq:2N1} and
Lemmas \ref{le:1} and \ref{le:2},
we see
\[
\begin{array}{@{}l}
L\left(\hatA_{3N_1},\hatB_{3N_1},\hatC_{3N_1},\hatD_{3N_1},
P_{3N_1},Q_{3N_1}, \gamma\right) 
\prec 0,\\
P_{3N_1}:=P_{N_1},\ 
Q_{3N_1}:=\left[
 \begin{array}{ccc}
  Q_{N_1} & 0 & 0 \\
  0 & Q_{N_1} & 0 \\
  0 & 0 & Q_{N_1}
 \end{array}
 \right].  
\end{array}
\]
By repeating this procedure $p-1$ times, we can conclude that
\rec{eq:N2} holds with
\[
 P_{N_2}=P_{N_1},\ 
 Q_{N_2}=\diag\underbrace{(Q_1,\cdots,Q_1)}_{\mbox{$p$ times}}.  
\]
This completes the proof.  
\end{proofof}

\section*{Lemma in the Proof of \rth{th:SG_RNN}}

In the proof of \rth{th:SG_RNN}
we use the next lemma.  
\begin{lemma}\normalfont
For given 
$\tilz\in\bbR^{\nz}$, 
$\tilw\in\bbR^{\nw}$, 
$\tils\in\bbR^{\ns}$, 
$\tilv\in\bbR^{\nv}$ and
$a,b,c,d\in\bbR$, suppose
\begin{equation}
\|\tilz\|_2\le a \|\tils\|_2+b\|\tilv\|_2,\ 
\label{eq:z} 
\end{equation}
\begin{equation}
\|\tilw\|_2\le c \|\tils\|_2+d\|\tilv\|_2.  
\label{eq:w} 
\end{equation}
Then, we have
\begin{equation}
\left\|
 \left[
 \begin{array}{c}
  \tilz\\
  \tilw
 \end{array}
 \right]\right\|_2 \le \sqrt{2}
\left\|
 \left[
 \begin{array}{cc}
  a & b \\
  c & d 
 \end{array}
 \right]\right\|_2
\left\|
 \left[
 \begin{array}{c}
  \tils\\
  \tilv
 \end{array}
 \right]\right\|_2.  
\label{eq:ineq}
\end{equation}
\label{le:ineq}
\end{lemma}
\begin{proofof}{\rle{le:ineq}}
From \rec{eq:z}, we first obtain
\[
\scalebox{0.8}{$
\arraycolsep=0.5mm
\begin{array}{@{}lcl}
 \|\tilz\|_2 \le 
 \left[
 \begin{array}{cc}
  a & b\\
 \end{array}\right]
 \left[
 \begin{array}{c}
  \|\tils\|_2 \\
  \|\tilv\|_2 \\
 \end{array}\right]&\le&
\left\|
 \left[
 \begin{array}{cc}
  a & b\\
 \end{array}\right]\right\|_2
\left\|
 \left[
 \begin{array}{c}
  \|\tils\|_2 \\
  \|\tilv\|_2 \\
 \end{array}\right]\right\|_2
\le
\left\|
 \left[
 \begin{array}{cc}
  a & b\\
  c & d
 \end{array}\right]\right\|_2
\left\|
 \left[
 \begin{array}{c}
  \|\tils\|_2 \\
  \|\tilv\|_2 \\
 \end{array}\right]\right\|_2.  
\end{array}$}
\]
It follows that
\[
 \|\tilz\|_2^2 \le 
\left\|
 \left[
 \begin{array}{cc}
  a & b\\
  c & d
 \end{array}\right]\right\|_2^2
\left\|
 \left[
 \begin{array}{c}
  \|\tils\|_2 \\
  \|\tilv\|_2 \\
 \end{array}\right]\right\|_2^2
 =
\left\|
 \left[
 \begin{array}{cc}
  a & b\\
  c & d
 \end{array}\right]\right\|_2^2
\left\|
 \left[
 \begin{array}{c}
  \tils \\
  \tilv \\
 \end{array}\right]\right\|_2^2.  
\]
Similarly, we have from \rec{eq:w} that
\[
 \|\tilw\|_2^2 \le 
\left\|
 \left[
 \begin{array}{cc}
  a & b\\
  c & d
 \end{array}\right]\right\|_2^2
\left\|
 \left[
 \begin{array}{c}
  \tils \\
  \tilv \\
 \end{array}\right]\right\|_2^2.  
\]
Therefore we have
\[
\left\|
 \left[
 \begin{array}{c}
  \tilz \\
  \tilw \\
 \end{array}\right]\right\|_2^2  
 =  \|\tilz\|_2^2+\|\tilw\|_2^2 
 \le 
2\left\|
 \left[
 \begin{array}{cc}
  a & b\\
  c & d
 \end{array}\right]\right\|_2^2
\left\|
 \left[
 \begin{array}{c}
  \tils \\
  \tilv \\
 \end{array}\right]\right\|_2^2.  
\]
This clearly shows that \rec{eq:ineq} holds.  
\end{proofof}

\end{document}